\documentclass[final,3p]{elsarticle}
\usepackage{mathrsfs,enumerate,amssymb,amsmath,color,lineno}
\usepackage[all]{xy}
\usepackage{latexsym}
\usepackage{amsthm,a4wide,color}
\usepackage{amscd,verbatim}
\usepackage{graphicx}
\usepackage[colorlinks=true]{hyperref}

\newtheorem{theorem}{Theorem}

\newtheorem{example}[theorem]{Example}
\newtheorem{question}[theorem]{Question}

\newtheorem{remark}[theorem]{Remark}

\journal{Topology and its Applications}

\begin{document}

\begin{frontmatter}

\title{A note on the sensitivity of semiflows\tnoteref{mytitlenote}}
\tnotetext[mytitlenote]{This work was supported by the National Natural Science Foundation of China
(No. 11601449), the Hong Kong Scholars Program, the National Natural Science Foundation of China (Key Program)
(No. 51534006), Science and Technology Innovation Team of Education Department of Sichuan for Dynamical System and its
Applications (No. 18TD0013), and Youth Science and Technology Innovation Team of Southwest
Petroleum University for Nonlinear Systems (No. 2017CXTD02).}

\author[a1]{Xinxing Wu\corref{mycorrespondingauthor}}
\cortext[mycorrespondingauthor]{Corresponding author}
\address[a1]{School of Sciences, Southwest Petroleum University, Chengdu, Sichuan 610500, China}
\ead{wuxinxing5201314@163.com}

\author[a2]{Xin Ma}
\address[a2]{bSchool of Science, Southwest University of Science and Technology, Mianyang, Sichuan 621010, China}
\ead{cauchy7203@gmail.com}

%
\author[a3]{Guanrong Chen}
\address[a3]{Department of Electronic Engineering, City University of Hong Kong,
Hong Kong SAR, China}
\ead{gchen@ee.cityu.edu.hk}

\author[a4]{Tianxiu Lu}
\address[a4]{School of Mathematics and Statistics, Sichuan University of Science and Engineering, Zigong, Sichuan 643000, China}
\ead{lubeeltx@163.com}

\begin{abstract}
In this note, it is shown that there exist two non-syndetically sensitive cascades defined on complete metric spaces whose product is cofinitely sensitive,
answering negatively the Question 9.2 posed in \cite[Miller, A., Money, C., Turk. J. Math., 41 (2017): 1323--1336]{MM2017}. Moreover, it is
shown that there exists a syndetically sensitive semiflow $(G, X)$ defined on a complete metric space $X$ such that
$(G_1, X)$ is not sensitive for some syndetic closed submonoid $G_1$ of $G$, answering negatively the Open question 3 posed
in \cite[Money, C., PhD thesis, University of Louisville, 2015]{M2015} and Question~43 posed in \cite[Miller, A., Real Anal. Exchange, 42 (2017): 9--24]{M2017}.
\end{abstract}

\end{frontmatter}

Let $\mathbb{N}=\left\{1, 2, 3, \ldots\right\}$, $\mathbb{N}_0=\left\{0, 1, 2, \ldots\right\}$ and
$G$ denote a noncompact abelian topological monoid with the identity element $0$. A subset $A$ of $G$ is
\begin{enumerate}[(1)]
\item {\it syndetic} if there exists a compact subset $K$ of $G$ such that, for every $t\in G$, $(t+K)\cap A\neq \O$,
where $t+K=\left\{t+a: a\in K\right\}$; 
\item {\it thick} if, for every compact subset $K$ of $G$, there exists some $t\in G$ such that $t+K\subset A$.
\end{enumerate}
Clearly, a subset $A$ of $G$ is syndetic (thick) if and only if $G\setminus A$ is not thick (not syndetic).

A jointly continuous monoid action $\pi: G\times X\longrightarrow X$ of $G$ on a metric space $(X, d)$ is called a {\it semiflow}
and denoted by $(G, X, \pi)$ or $(G, X)$. In particular, if the acting topological monoid is a topological group,
a semiflow is called a {\it flow}. The element $\pi (t, x)$ will be denoted by $t. x$ or $tx$, so that the defining
conditions for a semiflow have the following form:
$$
s.(t.x)=(s+t).x \text{ and } 0.x=x, \text{ for all } s, t\in G \text{ and all } x\in X.
$$
The maps
\begin{align*}
\pi_{t}: X &\longrightarrow X\\
x &\longmapsto t.x
\end{align*}
are called {\it transition maps}. For any $x\in X$, the set $Gx:=\{t.x: t\in G\}$ is called the {\it orbit} of $x$.
A semiflow $(G, X)$ is {\it minimal} if the orbit of every point $x\in X$ is dense in $X$, i.e., $\overline{Gx}=X$.
Otherwise, it is called {\it non-minimal}. For any subset $A$ of $X$ and any $t\in G$, let $t.A=\{t.x : x\in A\}$.
Now, let $X$ be a metric space and let $f: X\longrightarrow X$ be a continuous map, which leads to a natural semiflow,
where $G=\mathbb{N}_0$ (with the discrete topology), and for any $x\in X$ and any $n\in \mathbb{N}_0$, $n. x=f^{n}(x)$.
This type of semiflow is called a {\it cascade} and is often denoted by $(X, f)$ instead of $(\mathbb{N}_0, X)$.

Let $(G, X)$ be a semiflow. For any $\varepsilon>0$ and any subset $A$ of $X$, let
$$
\mathrm{N}(A, \varepsilon)=\left\{t\in G: \mathrm{diam}(t.A)>\varepsilon\right\}, \text{ where }
\mathrm{diam}(t.A) \text{ denotes the diameter of } t.A.
$$

A semiflow $(G, X)$ is
\begin{enumerate}[(1)]
\item {\it sensitive} if there exists $\varepsilon>0$ such that, for every nonempty open subset $U$
of $X$, $\mathrm{N}(U, \varepsilon)\neq \O$;
\item {\it syndetically sensitive} if there exists $\varepsilon>0$ such that, for every nonempty open subset $U$
of $X$, $\mathrm{N}(U, \varepsilon)$ is syndetic;
\item {\it thickly sensitive} if there exists $\varepsilon>0$ such that, for every nonempty open subset $U$
of $X$, $\mathrm{N}(U, \varepsilon)$ is thick;
\item {\it cofinitely sensitive} if there exists $\varepsilon>0$ such that, for every nonempty open subset $U$
of $X$, $\mathrm{N}(U, \varepsilon)$ is cofinite;
\item {\it multi-sensitive} if there exists $\varepsilon>0$ such that, for any nonempty open subsets $U_1, \ldots, U_k$
of $X$, $\bigcap_{i=1}^{k}\mathrm{N}(U_i, \varepsilon)\neq \O$.
\end{enumerate}

Recently, Miller and Money \cite{MM2015} proved that a non-minimal syndetically transitive semiflow is syndetically sensitive, generalizing some results
in \cite{BBCDS1992,GW1993,M2007,WLF2012}. Then, they \cite{MM2017} generalized some results on chaotic properties of cascades to the product of semiflows
and asked the following question. (For more recent results on sensitivity, refer to \cite{DK2010,KM2008,LOW2017,M2011,M2017,M2018-2,M2018-1,WC2017,WLML2019,WOC2016,WZ2012,WZ2013}
and some references therein.)
%
%

\begin{question}\label{Q-1}{\rm \cite[Question~9.2]{MM2017}
Let $(G, X)$ and $(G, Y)$ be two semiflows defined on two metric spaces $X$ and $Y$.
If $(G, X\times Y)$ is syndetically sensitive (resp., cofinitely sensitive,
multi-sensitive), is one of the factors syndetically sensitive (resp., cofinitely sensitive,
multi-sensitive)?}
\end{question}

We proved \cite{WWC2015} that $(X\times Y, f\times g)$ is multi-sensitive if and only if
$(X, f)$ or $(Y, g)$ is multi-sensitive. Using the analogous method in \cite{WWC2015}, Thakur and Das~\cite{TD2019}
showed that this also holds for semiflows. Money \cite{M2015} then conducted a systemic investigation
on chaotic properties for semiflows and asked the following question which also can be found in \cite[Question~43]{M2017}:

\begin{question}\label{Q-2}{\rm \cite[Open question~3]{M2015}
Let $(G, X)$ be a sensitive semiflow and $G_1$ be a syndetic closed submonoid of $G$.
Is $(G_1, X)$ necessarily sensitive?}
\end{question}

{\color{blue}Recently, Miller \cite{M2018-1} discussed relations between various types of sensitivity in general
semiflows and asked the following question:
\begin{question}\label{Q-3}{\rm \cite[Question~3]{M2018-1}
Does thick sensitivity imply syndetic sensitivity?}
\end{question}
}

This paper constructs two examples (Examples~\ref{Exa-1} and \ref{Exa-2} below),
answering negatively Questions \ref{Q-1}, \ref{Q-2}, and \ref{Q-3} above.

\begin{example}\label{Exa-1}
{\rm Let $L_0=\mathscr{L}_0=0$, $L_1=\mathscr{L}_1=2$, and $L_n=2^{L_1+ \cdots +L_{n-1}}\cdot (2n)$,
$\mathscr{L}_{n}=L_1+L_2+ \cdots +L_n$ for all $n\geq 2$, and set
$$
\begin{aligned}
X=[0, 1]&\bigcup \left(\bigcup_{n=1}^{+\infty}\bigcup_{i=0}^{2^{\mathscr{L}_{2n-2}}}\left[\mathscr{L}_{2n-1}+i,
\mathscr{L}_{2n-1}+i+\frac{1}{2n}\right]\right)\\
&\bigcup \left(\bigcup_{n=1}^{+\infty}\bigcup_{i=0}^{2^{\mathscr{L}_{2n-1}}}
\left[\mathscr{L}_{2n}+(2n+1)i, \mathscr{L}_{2n}+(2n+1)i+2n\right]\right),
\end{aligned}
$$
and
$$
\begin{aligned}
Y=[0, 1]&\bigcup \left(\bigcup_{n=1}^{+\infty}\bigcup_{i=0}^{2^{\mathscr{L}_{2n-2}}+4n-4}\left[\mathscr{L}_{2n-1}+2ni,
\mathscr{L}_{2n-1}+2ni+2n-1\right]\right)\\
&\bigcup \left(\bigcup_{n=1}^{+\infty}\bigcup_{i=0}^{2^{\mathscr{L}_{2n-1}}-4n+2}
\left[\mathscr{L}_{2n}+i, \mathscr{L}_{2n}+i+\frac{1}{2n}\right]\right).
\end{aligned}
$$

For $n\in \mathbb{N}$, let
$$
\mathscr{A}_n=\left[\mathscr{L}_{2n-1}+2^{\mathscr{L}_{2n-2}},
\mathscr{L}_{2n-1}+2^{\mathscr{L}_{2n-2}}+\frac{1}{2n-1}\right],
$$
$$
\mathscr{B}_n=\left[\mathscr{L}_{2n}+(2n+1)\cdot 2^{\mathscr{L}_{2n-1}}, \mathscr{L}_{2n}+(2n+1)\cdot 2^{\mathscr{L}_{2n-1}}+2n\right],
$$
$$
\mathscr{C}_n=\left[\mathscr{L}_{2n-1}+2n\cdot \left(2^{\mathscr{L}_{2n-2}}+4n-4\right),
\mathscr{L}_{2n-1}+2n\cdot \left(2^{\mathscr{L}_{2n-2}}+4n-4\right)+(2n-1)\right],
$$
and
$$
\mathscr{D}_{n}=\left[\mathscr{L}_{2n}+2^{\mathscr{L}_{2n-1}}-4n+2, \mathscr{L}_{2n}+2^{\mathscr{L}_{2n-1}}-4n+2+\frac{1}{2n}\right].
$$
Define $f: X\longrightarrow X$ and $g: Y\longrightarrow Y$ respectively by
$$
f(x)=\left\{
\begin{aligned}
&\frac{1}{2}x+2, \text{ if }  x\in [0, 1], \\
&x+1, \text{ if }  x\in \left[\mathscr{L}_{2n-1}+i, \mathscr{L}_{2n-1}+i+\frac{1}{2n}\right]
\text{ for some } 0\leq i<2^{\mathscr{L}_{2n-2}}, n\in \mathbb{N}, \\
&x+2n+1, \text{ if }  x\in \left[\mathscr{L}_{2n}+(2n+1)i, \mathscr{L}_{2n}+(2n+1)i+2n\right]
\text{ for some } 0\leq i<2^{\mathscr{L}_{2n-1}}, n\in \mathbb{N}, \\
&2n(2n-1)\left(x-
\mathscr{L}_{2n-1}-2^{\mathscr{L}_{2n-2}}\right)+\mathscr{L}_{2n}, \text{ if } x\in \mathscr{A}_n, n\in \mathbb{N},\\
&\frac{1}{2n(2n+1)}\left(x-\mathscr{L}_{2n}-(2n+1)\cdot 2^{\mathscr{L}_{2n-1}}\right)+\mathscr{L}_{2n+1},
\text{ if } x\in \mathscr{B}_n, n\in \mathbb{N},
\end{aligned}
\right.
$$
and
$$
g(x)=\left\{
\begin{aligned}
&x+2, \text{ if }  x\in [0, 1], \\
&x+1, \text{ if }  x\in \left[\mathscr{L}_{2n}+i, \mathscr{L}_{2n}+i+\frac{1}{2n}\right]
\text{ for some } 0\leq i<2^{\mathscr{L}_{2n-1}}-2n, n\in \mathbb{N}, \\
&x+2n, \text{ if }  x\in \left[\mathscr{L}_{2n-1}+2ni, \mathscr{L}_{2n-1}+2ni+(2n-1)\right]
\text{ for some } 0\leq i<2^{\mathscr{L}_{2n-2}}+2n-2, n\in \mathbb{N}, \\
&\frac{1}{2n (2n-1)}\left(x-\mathscr{L}_{2n-1}-2n\cdot \left(2^{\mathscr{L}_{2n-2}}+4n-4\right)\right)+\mathscr{L}_{2n},
\text{ if } x\in \mathscr{C}_n, n\in \mathbb{N},\\
&2n (2n+1) \left(x-\mathscr{L}_{2n}-2^{\mathscr{L}_{2n-1}}+4n-2\right)+\mathscr{L}_{2n+1},
\text{ if } x\in \mathscr{D}_n, n\in \mathbb{N}.
\end{aligned}
\right.
$$
Clearly, $f$ and $g$ are continuous. Arrange all closed intervals $[0, 1], [\mathscr{L}_1, \mathscr{L}_1+\frac{1}{2}], \ldots,
[\mathscr{L}_{1}+2^{\mathscr{L}_0}, \mathscr{L}_{1}+2^{\mathscr{L}_0}+\frac{1}{2}], [\mathscr{L}_2, \mathscr{L}_{2}+2], \ldots,
[\mathscr{L}_2+2\cdot 2^{\mathscr{L}_1}, \mathscr{L}_2+2\cdot 2^{\mathscr{L}_1}+2], \ldots, [\mathscr{L}_{2n-1},
\mathscr{L}_{2n-1}+\frac{1}{2n}], \ldots, [\mathscr{L}_{2n-1}+2^{\mathscr{L}_{2n-2}},
\mathscr{L}_{2n-1}+2^{\mathscr{L}_{2n-2}}+\frac{1}{2n}], [\mathscr{L}_{2n}, \mathscr{L}_{2n}+2n], \ldots, [\mathscr{L}_{2n}+(2n+1)2^{\mathscr{L}_{2n-1}},
\mathscr{L}_{2n}+(2n+1)2^{\mathscr{L}_{2n-1}}+2n], \ldots$ of $X$ by this natural order and denote them as $I_0, I_1, I_2, \ldots$. It is easy to see that
$I_{2k+2+2^{\mathscr{L}_1}+\cdots +2^{\mathscr{L}_{2k-1}}}=[\mathscr{L}_{2k+1}, \mathscr{L}_{2k+1}+\frac{1}{2(k+1)}]$, $\ldots$,
$I_{2k+2+2^{\mathscr{L}_1}+\cdots +2^{\mathscr{L}_{2k}}}=[\mathscr{L}_{2k+1}+2^{\mathscr{L}_{2k}},
\mathscr{L}_{2k+1}+2^{\mathscr{L}_{2k}}+\frac{1}{2(k+1)}]$, $I_{2k+3+2^{\mathscr{L}_1}+\cdots +2^{\mathscr{L}_{2k}}}=
[\mathscr{L}_{2k+2}, \mathscr{L}_{2k+2}+2(k+1)]$, $\ldots$, $I_{2k+3+2^{\mathscr{L}_1}+\cdots +2^{\mathscr{L}_{2k+1}}}=
[\mathscr{L}_{2k+2}+(2k+3)\cdot 2^{\mathscr{L}_{2k+1}}, \mathscr{L}_{2k+2}+(2k+3)\cdot 2^{\mathscr{L}_{2k+1}}+2(k+1)]$.
Similarly, arrange all closed intervals of $Y$ by this natural order and denote them as $J_0, J_1, J_2, \ldots$. It is easy to see that
$J_{2+2^{\mathscr{L}_1}+\cdots +2^{\mathscr{L}_{2k-1}}}=[\mathscr{L}_{2k+1}, \mathscr{L}_{2k+1}+(2k+1)]$, $\cdots$,
$J_{2+4k+2^{\mathscr{L}_1}+\cdots +2^{\mathscr{L}_{2k}}}=[\mathscr{L}_{2k+1}+2(k+1)\cdot (2^{\mathscr{L}_{2k}}+4k), \mathscr{L}_{2k+1}+2(k+1)\cdot (2^{\mathscr{L}_{2k}}+4k)+(2k+1)]$, $J_{3+4k+2^{\mathscr{L}_1}+\cdots +2^{\mathscr{L}_{2k}}}=[\mathscr{L}_{2k+2}, \mathscr{L}_{2k+2}+\frac{1}{2(k+1)}]$,
$\cdots$, $J_{1+2^{\mathscr{L}_1}+\cdots +2^{\mathscr{L}_{2k+1}}}=[\mathscr{L}_{2k+2}+2^{\mathscr{L}_{2k+1}}-4k-2, \mathscr{L}_{2k+2}+2^{\mathscr{L}_{2k+1}}-4k-2+\frac{1}{2(k+1)}]$.
Note that $f$ is a linear homeomorphism from $I_{n}$ to $I_{n+1}$, and $g$ is also a linear homeomorphism from $J_{n}$ to $J_{n+1}$,
for all $n\in \mathbb{N}_0$. According to the constructions of $f$ and $g$, it can be verified that
\begin{enumerate}[(i)]

\item $f$ and $g$ are continuous;

\item\label{i} $X$ and $Y$ are complete metric subspaces of $\mathbb{R}^1$,
implying that $X\times Y$ is a complete subspace of $\mathbb{R}^2$;


\item\label{iv} for any $n\in \left[2^{\mathscr{L}_1}+\cdots + 2^{\mathscr{L}_{2k-1}}+2k+2, 2^{\mathscr{L}_1}+\cdots +
2^{\mathscr{L}_{2k}}+2k+2\right]$ ($k\in \mathbb{N}$),
$$
f^{n}([0, 1])=\left[\mathscr{L}_{2k+1}+n-a_k, \mathscr{L}_{2k+1}+n-a_k+\frac{1}{2k+2}\right],
$$
where $a_k=2^{\mathscr{L}_1}+\cdots + 2^{\mathscr{L}_{2k-1}}+2k+2$,
implying that $\mathrm{diam}(f^{n}([0,1]))=\frac{1}{2k+2}$;

\item\label{v} for any $n\in \left[2^{\mathscr{L}_1}+\cdots + 2^{\mathscr{L}_{2k}}+2k+3, 2^{\mathscr{L}_1}+\cdots +
2^{\mathscr{L}_{2k+1}}+2k+3\right]$ ($k\in \mathbb{N}$),
$$
f^{n}([0, 1])=\left[\mathscr{L}_{2k+2}+(n-b_k)(2k+3), \mathscr{L}_{2k+2}+(n-b_k)(2k+3)+2k+2\right],
$$
where $b_k=2^{\mathscr{L}_1}+\cdots +
2^{\mathscr{L}_{2k}}+2k+3$, implying that $\mathrm{diam}(f^{n}([0,1]))=2k+2$;

\item\label{vi} for any $n\in \left[2^{\mathscr{L}_1}+\cdots + 2^{\mathscr{L}_{2k-1}}+2, 2^{\mathscr{L}_1}+\cdots +
2^{\mathscr{L}_{2k}}+4k+2\right]$ ($k\in \mathbb{N}$),
$$
g^{n}([0, 1])=\left[\mathscr{L}_{2k+1}+(n-c_k)(2k+2), \mathscr{L}_{2k+1}+(n-c_k)(2k+2)+2k+1\right],
$$
where $c_k=2^{\mathscr{L}_1}+\cdots + 2^{\mathscr{L}_{2k-1}}+2$,
implying that $\mathrm{diam}(g^{n}([0,1]))=2k+1$;

\item\label{vii} for any $n\in \left[2^{\mathscr{L}_1}+\cdots + 2^{\mathscr{L}_{2k}}+4k+3, 2^{\mathscr{L}_1}+\cdots +
2^{\mathscr{L}_{2k+1}}+1\right]$ ($k\in \mathbb{N}$),
$$
g^{n}([0, 1])=\left[\mathscr{L}_{2k+2}+n-d_k, \mathscr{L}_{2k+2}+n-d_k+\frac{1}{2k+2}\right],
$$
where $b_k=2^{\mathscr{L}_1}+\cdots + 2^{\mathscr{L}_{4k}}+2k+3$,
implying that $\mathrm{diam}(g^{n}([0,1]))= \frac{1}{2k+2}$.
\end{enumerate}

\medskip

{\bf Claim 1.} $(X, f)$ is not syndetically sensitive.

\medskip

Fix an open subset $U=[0, 1)$ of $X$. For any $\varepsilon>0$, applying \eqref{iv}
yields that
$\{n\in \mathbb{N}_0: \mathrm{diam}(f^n(U))\leq \varepsilon\} \supset
\bigcup_{n=K}^{+\infty}\left[2^{\mathscr{L}_1}+\cdots + 2^{\mathscr{L}_{2k-1}}+2k+2, 2^{\mathscr{L}_1}+\cdots +
2^{\mathscr{L}_{2k}}+2k+2\right]$ for some $K\in \mathbb{N}$, i.e., $\{n\in \mathbb{N}_0: \mathrm{diam}(f^n(U))\leq \varepsilon\}$
is a thick set. This implies that $\mathrm{N}_{f}(U, \varepsilon)$ is not syndetic. Therefore, $(X, f)$ is not syndetically
sensitive.

\medskip

{\bf Claim 2.} $(Y, g)$ is not syndetically sensitive.

\medskip

Similarly to the proof of Claim 1, applying \eqref{vii} follows that this is true.

\medskip

{\bf Claim 3.}  $(X\times Y, f\times g)$ is cofinitely sensitive.

\medskip

Given any nonempty open subset $W$ of $X\times Y$, there exist open subsets $U \subset X$
and $V\subset Y$ such that $U\times V\subset W$. From \eqref{iv}--\eqref{vii}, it follows that
there exist non-degenerate closed intervals $[\alpha_1, \beta_1], [\alpha_2, \beta_2] \subset [0, 1]$
and $p, q\in \mathbb{N}_0$ such that $f^p([\alpha_1, \beta_1])\subset U$ and $g^q([\alpha_2, \beta_2])\subset V$.
As $f$ and $g$ are piecewise linear mappings, for any $n\in \mathbb{N}_0$, one has
$$
\mathrm{diam}(f^n([\alpha_1, \beta_1]))=(\beta_1-\alpha_1)\cdot \mathrm{diam}(f^n([0, 1])),
$$
and
$$
\mathrm{diam}(g^n([\alpha_2, \beta_2]))=(\beta_2-\alpha_2)\cdot \mathrm{diam}(g^n([0, 1])).
$$
This, together with \eqref{v} and \eqref{vi}, implies that
there exists $K>2(p+q)$ such that
$$
\mathscr{P}:=\bigcup_{k=K}^{+\infty}\left[2^{\mathscr{L}_1}+\cdots + 2^{\mathscr{L}_{2k}}+2k+3-p, 2^{\mathscr{L}_1}+\cdots +
2^{\mathscr{L}_{2k+1}}+2k+3-p\right]\subset \mathrm{N}(U, 1),
$$
and
$$
\mathscr{Q}:=\bigcup_{k=K}^{+\infty}\left[2^{\mathscr{L}_1}+\cdots + 2^{\mathscr{L}_{2k-1}}+2-q, 2^{\mathscr{L}_1}+\cdots +
2^{\mathscr{L}_{2k}}+4k+2-q\right]\subset \mathrm{N}(V, 1),
$$
implying that
$$
\mathrm{N}(W, 1)\supset \mathrm{N}(U\times V, 1)\supset \mathrm{N}(U, 1)
\cup \mathrm{N}(V, 1)\supset \mathscr{P}\cup \mathscr{Q}.
$$
Noting that $2^{\mathscr{L}_1}+\cdots +2^{\mathscr{L}_{2k}}+4k+2-q>2^{\mathscr{L}_1}+\cdots +2^{\mathscr{L}_{2k}}+2k+3-p$
and $2^{\mathscr{L}_1}+\cdots +2^{\mathscr{L}_{2k+1}}+2k+3-p>2^{\mathscr{L}_1}+\cdots +2^{\mathscr{L}_{2k+1}}+2-q$ hold for all $k\geq K$,
it can be verified that
$$
\left[2^{\mathscr{L}_1}+\cdots + 2^{\mathscr{L}_{2K-1}}+2-q, +\infty\right) = \mathscr{P}\cup \mathscr{Q},
$$
i.e., $\mathscr{P}\cup \mathscr{Q}$ is a cofinite set. This implies that $f\times g$ is cofinitely sensitive.

\medskip

{\color{blue}{\bf Claim 4.} $(X, f)$ is thickly sensitive.

\medskip

Given any nonempty open subset $U$ of $X$, similarly to the proof of Claim 3, one has
\begin{enumerate}[(i)]
\item there exist non-degenerate closed interval $[\alpha, \beta]\subset [0, 1]$
and $p\in \mathbb{N}_0$ such that $f^p([\alpha, \beta])\subset U$;
\item for any $n\geq p$,
$$
\mathrm{diam}(f^{n}(U))\geq \mathrm{diam}(f^{n}(f^p([\alpha, \beta])))
=\mathrm{diam}(f^{n+p}([\alpha, \beta]))=(\beta-\alpha)\cdot \mathrm{diam}(f^{n+p}([0, 1])).
$$
\end{enumerate}
This, together with \eqref{v}, implies that there exists $K\in \mathbb{N}$ such that
$$
\bigcup_{k=K}^{+\infty}\left[2^{\mathscr{L}_1}+\cdots + 2^{\mathscr{L}_{2k}}+2k+3-p, 2^{\mathscr{L}_1}+\cdots +
2^{\mathscr{L}_{2k+1}}+2k+3-p\right]\subset \mathrm{N}(U, 1),
$$
implying that $(X, f)$ is thickly sensitive.
}
}
\end{example}

\begin{remark}
{\rm
\begin{enumerate}[(1)]
\item Example \ref{Exa-1} shows that the answer to Question \ref{Q-1} for syndetic sensitivity and cofinite sensitivity
is negative. Combining this with \cite[Theorem~3.6, Proposition~3.7]{TD2019} completely solves Question \ref{Q-1}.
\item Example \ref{Exa-1}, together with \cite[Theorem~10]{WWC2015}, also gives a complete answer to \cite[Remark~3.4]{LZ2013}.
{\color{blue}\item From Claims 1 and 4 of Example \ref{Exa-1}, it follows that there exists a thickly sensitive semiflow which is not syndetically
sensitive, answering negatively Question~\ref{Q-3}. }
\end{enumerate}
}
\end{remark}

Similarly to the constructed semiflow in Example \ref{Exa-1}, the following example gives an negative answer to Question \ref{Q-2},
showing that there exists a syndetically sensitive semiflow $(\mathbb{N}_0, X)$ defined on a complete metric space $X$ such
that $(2\mathbb{N}_0, X)$ is not sensitive.

\begin{example}\label{Exa-2}
{\rm Let $L_0=\mathscr{L}_0=0$, $L_1=\mathscr{L}_1=2$, and $L_n=2^{L_1+ \cdots +L_{n-1}}$,
$\mathscr{L}_{n}=L_1+L_2+ \cdots +L_n$ for all $n\geq 2$ and set
$X=[0, 1]\cup \left(\bigcup_{n=1}^{+\infty}\left[\mathscr{L}_{2n-1}, \mathscr{L}_{2n-1}+2n-1\right]\right)\cup
\left(\bigcup_{n=1}^{+\infty}\left[\mathscr{L}_{2n}, \mathscr{L}_{2n}+\frac{1}{2n}\right]\right)$. Define
$f: X\longrightarrow X$ as
$$
f(x)=\left\{
\begin{aligned}
&x+2, \text{ if }  x\in [0, 1], \\
&\frac{1}{2n(2n-1)}\left(x-\mathscr{L}_{2n-1}\right)+\mathscr{L}_{2n}, \text{ if }  x\in \left[\mathscr{L}_{2n-1}, \mathscr{L}_{2n-1}+2n-1\right]
\text{ for some } n\in \mathbb{N}, \\
&2n(2n+1)\left(x-\mathscr{L}_{2n}\right)+\mathscr{L}_{2n+1}, \text{ if }  x\in \left[\mathscr{L}_{2n}, \mathscr{L}_{2n}+\frac{1}{2n}\right]
\text{ for some } n\in \mathbb{N}.
\end{aligned}
\right.
$$
For any nonempty open subset $U$ of $X$, there exist a non-degenerate closed interval $I \subset [0, 1]$ and $k\in \mathbb{N}_0$
such that $f^{k}(I)\subset U$. Thus, for any $n>k$, one has
$$
\mathrm{diam}(f^{2n+1-k}(U))\geq \mathrm{diam} (f^{2n+1}(I))=(2n+1)\cdot |I|
\longrightarrow + \infty\quad (n\longrightarrow +\infty).
$$
This implies that $(X, f)$ is syndetically sensitive.

Meanwhile, it is easy to see that, for any $n\in \mathbb{N}$, $\mathrm{diam}(f^{2n}([0, 1]))=\frac{1}{2n}$, implying that $(X, f^2)=(2\mathbb{N}_0, X)$
is not sensitive.
}
\end{example}


\section*{References}
\bibliographystyle{abbrv}
\bibliography{mybibfile}{}
\end{document}